\documentclass[11pt,reqno]{amsart}
\usepackage{enumerate}
\usepackage{url}
\usepackage{cite}
\usepackage{comment}
\usepackage{graphicx}
\usepackage{tikz}
\usetikzlibrary{arrows,shapes,trees,backgrounds}
\usepackage{amsfonts, amsmath, amssymb, amscd, amsthm, bm, cancel}
\usepackage[
linktocpage=true,colorlinks,citecolor=magenta,linkcolor=blue,urlcolor=magenta]{hyperref}
 \newtheorem{theorem}{Theorem}[section]
 \newtheorem{corollary}[theorem]{Corollary}
  
 \newtheorem{lemma}[theorem]{Lemma}

 \theoremstyle{definition}
 \newtheorem{definition}[theorem]{Definition}
  
 \theoremstyle{remark}

 \numberwithin{equation}{section}
 \allowdisplaybreaks

\begin{document}
\title[Subgradient estimates for a nonlinear subparablic equation]{Subgradient estimates for a nonlinear subparablic equation on complete pseudo-Hermitian manifolds}
\author[Wenjing Wu]{Wenjing Wu}
\address{School of Mathematical Sciences, University of Science and Technology of China, Hefei 230026, P.R. China}
\email{\href{mailto:anprin@mail.ustc.edu.cn}{anprin@mail.ustc.edu.cn}}

\date{\today}

\begin{abstract}
Let $(M,J,\theta)$ be a complete pseudo-Hermitian manifold which satisfies the CR sub-Laplacian comparison property. In this paper, we derive the local subgradient estimates for positive solutions to the following nonlinear subparabolic equation:
\begin{equation*}
u_t=\Delta_\mathrm{b}u+au\log u+bu,
\end{equation*}
on $M\times [0,+\infty)$, where $a,b\in \mathbb{R}$. As a application, we derive a priori estimate and a Harnack inequality for positive solutions to the subelliptic equation $\Delta_\mathrm{b}u+au\log u=0$.
\end{abstract}

\maketitle

\section{Introduction}
In \cite{LM}, L.Ma proved local gradient estimates of positive solutions to the equation
\begin{equation}\label{0.0}
\Delta u+au\log u+bu=0,
\end{equation}
for complete noncompact Riemannian manifolds with a fixed metric and Ricci curvature locally bounded below, where $a<0$ and $b\in \mathbb{R}$ are constants. In \cite{Y}, Yang generalized L.Ma's result and derived local gradient estimates for positive solutions to the parabolic equation
\begin{equation}\label{eq}
u_t=\Delta u+au\log u+bu,
\end{equation}
where $a,b\in \mathbb{R}$ are constants.

\begin{theorem}[\cite{Y}]\label{thmA}
Let $M$ be a complete noncompact Riemannian manifold of dimension $n$. Suppose the Ricci curvature of $M$ is bounded below by the constant $-K$, where $R>0,K \geq 0$, in the geodesic ball $B_{2R}(p)$ with radius $2R$ around $p \in M$. If $u(x, t)$ is a positive smooth solution to equation $(\ref{0.0})$ on $M \times[0,+\infty)$, let $f=\log u$. Then\\
$(i)$ if $a<0$, we have for any $\alpha>1$ and $0<\delta<1$,
$$
\begin{gathered}
|\nabla f|^2(x, t)+\alpha a f(x, t)+\alpha b-\alpha f_t(x, t) \leq \frac{n \alpha^2}{2 \delta t}+\frac{n \alpha^2}{2 \delta}\left\{\frac{2 \epsilon^2}{R^2}+\frac{\nu}{R^2}-\frac{a}{2}\right. \\
\left.+\frac{\epsilon^2}{R^2}(n-1)(1+R \sqrt{K})+\frac{K}{\alpha-1}+\frac{n \alpha^2 \epsilon^2}{8(1-\delta)(\alpha-1) R^2}\right\},
\end{gathered}
$$
on $B_R(p) \times(0,+\infty)$, where $\epsilon>0$ and $\nu>0$ are some constants;\\
$(ii)$ if $a>0$, we have for any $\alpha>1$ and $0<\delta<1$,
$$
\begin{gathered}
|\nabla f|^2(x, t)+\alpha a f(x, t)+\alpha b-\alpha f_t(x, t) \leq \frac{n \alpha^2}{2 \delta t}+\frac{n \alpha^2}{2 \delta}\left\{\frac{2 \epsilon^2}{R^2}+\frac{\nu}{R^2}+a\right. \\
\left.+\frac{\epsilon^2}{R^2}(n-1)(1+R \sqrt{K})+\frac{K}{\alpha-1}+\frac{n \alpha^2 \epsilon^2}{8(1-\delta)(\alpha-1) R^2}\right\},
\end{gathered}
$$
on $B_R(p) \times(0,+\infty)$, where $\epsilon>0$ and $\nu>0$ are some constants.
\end{theorem}
In this paper, we consider the following nonlinear subparabolic equation
\begin{equation}\label{1.1}
u_t=\Delta_\mathrm{b}u+au\log u+bu,
\end{equation}
on a complete pseudo-Hermitian manifold $(M^{2n+1},J,\theta)$, where $a,b$ are two real constants. Replacing $u$ by $e^{b/a}u$, we only need to consider positive smooth solutions of the following equation:
\begin{equation}\label{equation}
u_t=\Delta_\mathrm{b}u+au\log u.
\end{equation}

To seek an appropriate method to deal with the locally positive solutions for the $(\ref{equation})$, let us recall what is the proof of Theorem $\ref{thmA}$. There are two main ingredients: the Bochner formula and a Laplacian comparison theorem. The Bochner formula states that 
$$\frac{1}{2}\Delta |\nabla \varphi|^2=|\nabla^2 \varphi|^2+\langle \nabla \varphi,\nabla \Delta \varphi\rangle +\operatorname{Ric}(\nabla \varphi,\nabla \varphi)$$
for any smooth function $\varphi$ on $M$. The Laplacian comparison theorem states that if $\operatorname{Ric}_M \geq -K$, then we have $ \Delta r\leq (n-1)r^{-1}+\sqrt{K(n-1)},$ in the distribution sense. By using the Bochner formula to $\log u$, one deduces a differential inequality 
\begin{align*}
	(\Delta-\partial_t) F
	\geq & -2 \langle \nabla f, \nabla F\rangle +t[2/n(-F/(\alpha t)-(1-\alpha^{-1})|\nabla f|^2)^2 \\
	& +((\alpha-1) a-2 K)|\nabla f|^2]-a F-\frac{F}{t},
\end{align*}
where $f=\log u$ and $F=t(|\nabla f|^2+\alpha a f-\alpha f_t)$. One obtain the local gradient estimates by combining the above inequality along with a localiztion argument by using a suitable cut-off function. The estimates on the cut-off function makes use of the Laplacian comparison theorem.

In this paper, we need to consider the analogues in CR world. Firstly, let us consider the Bochner formula in CR context. Greenleaf \cite{G} proved that the following CR version of Bochner formula in a pseudo-Hermitian $(2n+1)$-manifold: for a smooth real-valued function $\varphi$ on $M$,
\begin{align*}
\Delta_\mathrm{b} |\nabla^\mathrm{b} \varphi|^2 &= 2|\pi_\mathrm{b}\nabla^2 \varphi|^2+2\langle\nabla^\mathrm{b} \varphi,\nabla^\mathrm{b} \Delta_\mathrm{b}\varphi\rangle \\
& \quad+2(2\operatorname{Ric}-(n-2)\operatorname{Tor})(\nabla^{1,0}\varphi,\nabla^{1,0}\varphi)+4\langle J\nabla^\mathrm{b}\varphi,\nabla^\mathrm{b} \varphi_0\rangle.
\end{align*}
Notice that the right-hand side of the CR Bochner formula involves a term $\langle J\nabla^\mathrm{b}\varphi,\nabla^\mathrm{b} \varphi_0\rangle$ that has no analogue in the Riemannian case. In order to overcome this difficulty, we consider an auxiliary function $F=t(|\nabla f|^2+\alpha a f-\alpha f_t+\beta tf_0^2)$.

Next, we consider the sub-Laplacian comparison theorem in CR context. When a complete pseudo-Hermitian manifold of vanishing pseudo-Hermitian torsion satisfies 
$\operatorname{Ric}(Z,Z)\geq-k|Z|^2,\text{ for all }Z\in T_{1,0}(M),$ and $k\geq0$, Chang et al. \cite{CKLT} proved that the CR sub-Laplacian comparison theorem of the Carnot-Carath\'{e}odory distance function holds. On the other hand, Dong et al. \cite{CDRZ} proved the following comparison theorem:
\begin{theorem}[\cite{CDRZ}]
Suppose $(M^{2n+1},J,\theta)$ is a complete pseudo-Hermitian manifold. If for some $k_1,k_2\geq0$,
\begin{equation*}
R_*\geq-k_1,\text{ and }|A|,|div A|\leq k_2,\text{ on }M,
\end{equation*}
then there exists $C_1=C_1(n)$ such that
\begin{equation*}
\Delta_\mathrm{b} r\leq C_1\bigg(\frac{1}{r}+\sqrt{1+k_1+k_2+k_2^2}\bigg),
\end{equation*}
in the distribution sense, where $r(x)$ is the Riemannian distance from $x_0$. 
\end{theorem}
In particular, if $(M,J,\theta)$ has zero torsion and satisfies $\operatorname{Ric}(Z,Z)\geq-k|Z|^2$, for all $Z\in T_{1,0}(M),k\geq0($this implies $R_*\geq-k)$, the CR sub-Laplacian comparison theorem of the Riemannian distance function holds.

\vspace{0.5cm}

In CR geometry, Chang et al. \cite{CKLT} derived subgradient estimates for positive pseudoharmonic functions in a complete noncompact pseudo-Hermitian manifold which satisfies the CR sub-Laplacian comparison property. In \cite{HZ}, He and Zhao obtain subgradient estimates for positive solutions to the following nonlinear equation:
\begin{equation*}
\Delta_b u+cu^{-\alpha}=0,
\end{equation*}
where $c\in \mathbb{R}$ and $\alpha>-1$ are two constants. In \cite{HJL}, Y. B. Han, K. G. Jiang, and M. H. Liang modified the arguments of \cite{CKLT} and obtained the following result.
\begin{theorem}[\cite{HJL}]
Let $(M, J, \theta)$ be a complete pseudo-Hermitian $(2n+1)$-manifold. Suppose that
$$2 \operatorname{Ric}(Z,Z)-(n-2) \operatorname{Tor}(Z,Z) \geq-2 k|Z|^2,$$
for all $Z \in T_{1,0} (M)$ and $k\geq 0$. Furthermore, we assume that $(M,J,\theta)$ satisfies the CR sub-Laplacian comparison property. If $u$ is the positive solution of
$$
\Delta_\mathrm{b} u+a u \log u=0
$$
with $[\Delta_\mathrm{b}, T]u=0$ on $M$, let $f=\log u$. Then we have
\begin{align*}
& \frac{|\nabla_b u|^2}{u^2}+a \frac{n+a \beta+4+\gamma+2 \beta k}{n+a \beta} \log u+\beta \frac{u_0^2}{u^2} \\
< & \frac{(n+a \beta+4+\gamma+2 \beta k)^2}{2(4+\gamma+2 \beta k)^2}\left[2 k+\frac{4}{\beta}-a \frac{4+\gamma+2 \beta k}{n+a \beta}+\frac{C}{R}\right]
\end{align*}
on the ball $B_R(x_0)$ of large enough radius $R$ which depends only on $a,\beta, \gamma, k$, where $a \leq 0, \beta \geq 0$ and $\gamma>0$ are constants such that $n+a \beta>0$.	
\end{theorem}
	
In this paper, let $B_R(x_0)$ be the Riemannian geodesic ball of radius $R$ centered at $x_0\in M$.
\begin{theorem}\label{thm}
Let $(M^{2n+1},J,\theta)$ be a complete pseudo-Hermitian manifold with
\begin{equation*}
(2\operatorname{Ric}-(n-2)\operatorname{Tor})(Z,Z)\geq-2k|Z|^2,\text{ for all }Z\in T_{1,0}(M), \text{ on }B_{2R}(x_0),
\end{equation*}
where $k\geq0$. Suppose that $M$ satisfies the CR sub-Laplacian comparison property $(\ref{CRL})$.

Let $u(x,t)$ be a positive smooth solution to equation $(\ref{equation})$ on $B_{2R}(x_0)\times [0,+\infty)$, satisfying  
\begin{equation*}
[\Delta_\mathrm{b},T]u=0.
\end{equation*} 
Then\\
$(i)$ If $a\leq0$, we have for any $\alpha\geq4/n+1$ and $0<\delta<1$,
\begin{align}
\frac{|\nabla_b u|^2}{u^2}&+\alpha a f-\alpha\frac{u_t}{u}\leq
\frac{n\alpha^2}{\delta t}\bigg(1+\frac{2}{n(\alpha-1)}\bigg)\notag\\
&\quad+ \frac{n\alpha^2}{\delta}\bigg[\bigg(\frac{\upsilon}{R^2}+\frac{\epsilon C_0}{R^2}(1+\sqrt{\tilde{k}}R)+\frac{2\epsilon^2}{R^2}\bigg)+\frac{((\alpha-1)|a|+2k)}{\alpha-1}\notag\\
&\quad+\frac{\epsilon^2 \alpha^2 n}{2(1-\delta)(\alpha-1)R^2}\bigg].\label{0.1}
\end{align}
on $B_R(x_0)\times (0,+\infty)$, where $f=\log u$ and $\epsilon,\upsilon,C_0,\tilde{k}>0$ are some constants;\\
$(ii)$ If $a>0$, we have for any $\alpha\geq4/n+1$ and $0<\delta<1$,
\begin{align}
\frac{|\nabla_b u|^2}{u^2}&+\alpha a f-\alpha \frac{u_t}{u}\leq
\frac{ n\alpha^2}{\delta t}\bigg(1+\frac{2}{n(\alpha-1)}\bigg)\notag\\
&\quad+\frac{ n\alpha^2}{\delta }\bigg[\bigg(\frac{\upsilon}{R^2}+\frac{\epsilon C_0}{R^2}(1+\sqrt{\tilde{k}}R)+\frac{2\epsilon^2}{R^2}\bigg)+a+\frac{2kn+2a}{n(\alpha-1)}\notag\\
&\quad+\frac{\epsilon^2 \alpha^2 n}{2(1-\delta)(\alpha-1)R^2}\bigg],\label{0.2}
\end{align}
on $B_R(x_0)\times (0,+\infty)$, where $f=\log u$ and $\epsilon,\upsilon,C_0,\tilde{k}>0$ are some constants.
\end{theorem}

An interesting corollary of Theorem $\ref{thm}$ is the following result.
\begin{corollary}\label{cor}
Let $(M^{2n+1},J,\theta)$ be a complete pseudo-Hermitian manifold with
\begin{equation*}
(2Ric-(n-2)Tor)(Z,Z)\geq 0,\text{ for all }Z\in T_{1,0}(M), \text{ on }M.
\end{equation*}
Suppose that $M$ satisfies the CR sub-Laplacian comparison property $(\ref{CRL})$. If $u(x)$ is a positive smooth solution to the equation
\begin{equation*}
\Delta_\mathrm{b}u+au\log u=0,\;\text{ on }M,
\end{equation*}
satisfying 
\begin{equation*}
[\Delta_\mathrm{b},T]u=0.
\end{equation*} 
Then if $a<0$, we have $u(x)\geq e^{-4-n}$ for all $x\in M$; and if $a>0$, we have $u(x)\leq e^{3(4+n)/2}$ for all $x\in M$.
\end{corollary}
\begin{theorem}\label{thm2}
Under the assumptions of Theorem $\ref{thm}$, if $u$ is a positive solution to $(\ref{equation})$, then for every $(x_1,t_1),(x_2,t_2)$ in $B_{R/2}(x_0)\times (0,\tau)$ with $t_2>t_1$ and $\alpha \geq4/n+1$, we have
\begin{equation*}
u(x_2,t_2)\geq u(x_1,t_1)\bigg(\frac{t_2}{t_1}\bigg)^{-\frac{\alpha(n\alpha-n+2)}{\delta(\alpha-1) }}e^{-\alpha L(x_1,x_2,t_2-t_1)}e^{S(t_2-t_1)}
\end{equation*}
Here $S$ is a constant given by $(\ref{5.2})$. Furthermore $L$ is given by 
$$ L(x_1,x_2,t_2-t_1)=\inf_{\gamma}\bigg[\frac{1}{4(t_2-t_1)}\int_0^1 |\dot{\gamma}(t)|^2 dt\bigg],$$
where $\Gamma$ is the set of all lengthy curves $\gamma\in C^1([t_1,t_2];M)$ lying entirely in $B_{R}(x_0)$ with $\gamma(t_1)=x_1$ and $\gamma(t_2)=x_2$.
\end{theorem}
\section{Preliminaries}
In this section, we introduce some basic notions of pseudo-Hermitian geometry. For details, readers may refer to \cite{DT}.

Recall that a smooth manifold $M$ of real dimension $2n+1$ is said to be a CR manifold if there exists a smooth rank $n$ complex subbundle $T_{1,0}(M)\subset T^{\mathbb{C}}M$ such that
\begin{equation*}
T_{1,0}(M)\cap T_{0,1}(M)=\{0\}
\end{equation*}
\begin{equation*}
[\Gamma^\infty(T_{1,0}(M)),\Gamma^\infty(T_{1,0}(M))]\subset \Gamma^\infty(T_{1,0}(M)),
\end{equation*}
where $T_{0,1}(M)=\overline{T_{1,0}(M)}$ is the complex conjugate of $T_{1,0}(M)$. Equivalently, the CR structure may also be described by the real subbundle $H(M)=Re\{T_{1,0}(M)\oplus T_{0,1}(M)\}$ of $TM$ which carries an almost complex structure $J:H(M)\rightarrow H(M)$ defined by $J(X+\bar{X})=i(X-\bar{X})$ for any $X\in T_{1,0}(M)$. Since $H(M)$ is naturally oriented by the almost complex structure $J$, then $M$ is orientable if and only if there exists a global nowhere vanishing $1$-form $\theta$ such that $H(M)=Ker(\theta)$. Any such section $\theta$ is referred to as a pseudo-Hermitian structure on $M$. The Levi form $L_\theta$ of a given pseudo-Hermitian structure $\theta$ is defined by
\begin{equation*}
L_\theta(Z,W)=-id\theta(Z,\bar{W}) \text{ for any }Z,W\in T_{1,0}(M).
\end{equation*}
An orientable CR manifold $(M,J,\theta)$ is called strictly pseudo-convex if $L_\theta$ is positive definite for some $\theta$. Such a triple $(M,J,\theta)$ is called a pseudo-Hermitian manifold.

For a pseudo-Hermitian manifold $(M,J,\theta)$, there exists a unique nowhere zero vector field $T$, called the Reeb vector field, satisfying $\theta(T)=1,T\lrcorner d\theta=0$. It gives a decomposition of the tangent bundle $TM$:
\begin{equation*}
TM=\mathbb{R}T\oplus H(M),
\end{equation*}
which induces the projection $\pi_H:TM\rightarrow H(M)$. Define the bilinear form $G_\theta$ by setting
\begin{equation*}
G_\theta(X,Y)=d\theta(X,JY),
\end{equation*}
for any $X,Y\in H(M)$. Now, it is natural to define a Riemannian metric
\begin{equation*}
g_\theta=\pi_H G_\theta+\theta\otimes \theta
\end{equation*}
which makes $H(M)$ and $\mathbb{R}T$ orthogonal. The metric $g_\theta$ is called Webster metric, which is also denoted by $\langle\cdot,\cdot\rangle$ for simplicity. By requiring $JT=0$, the almost complex structure $J$ can be extended to an endomorphism of $TM$. In this paper, a pseudo-Hermitian manifold $(M,J,\theta)$ is called complete about the Webster metric $g_\theta$.

On a pseudo-Hermitian manifold, there exists a canonical connection $\nabla$, which is called Tanaka-Webster connection, preserving the horizontal distribution, almost complete structure and Webster metric.

We first define Ric and Tor on $T_{1,0}(M)$ by
\begin{equation*}
\operatorname{Ric}(Z,W)=R_{\alpha\bar{\beta}}Z^\alpha W^{\bar{\beta}}
\end{equation*}
and
\begin{equation*}
\operatorname{Tor}(Z,W)=i\sum_{\alpha,\beta}(A_{\bar{\alpha}\bar{\beta}}Z^{\bar{\alpha}} W^{\bar{\beta}}-A_{\alpha\beta}Z^\alpha Z^\beta).
\end{equation*}
Here $Z=Z^\alpha T_\alpha,W=W^\beta T_\beta$ for a frame $\{T,T_\alpha,T_{\bar{\alpha}}\}$ of $T^{\mathbb{C}}M$. $R_{\gamma \alpha\bar{\beta}}^\delta$ is the pseudo-Hermitian curvature tensor; $R_{\alpha\bar{\beta}}=R_{\gamma \alpha\bar{\beta}}^\gamma$ is the pseudo-Hermitian Ricci curvature tensor; and $A_{\alpha\beta}$ is the torsion tensor. Let $\{\theta,\theta^\alpha,\theta^{\bar{\alpha}}\}$ be the dual frame of $\{T,T_\alpha,T_{\bar{\alpha}}\}$. Then, we can write the Levi form by
\begin{equation*}
L_\theta=h_{\alpha\bar{\beta}}\theta^\alpha\wedge \theta^{\bar{\beta}}
\end{equation*}
for some positive definite hermitian matrix of functions $(h_{\alpha\bar{\beta}})$. Actually, we can always choose $T_\alpha$ such that $h_{\alpha\bar{\beta}}=\delta_{\alpha\beta}$; hence, throughout this paper, we assume $h_{\alpha\bar{\beta}}=\delta_{\alpha\beta}$.

Since the Tanaka-Webster connection parallelizes the $T_{1,0}(M)$ there exists uniquely defined complex $1$-forms $\omega_\beta^\alpha\in \Gamma(T^{*\mathbb{C}}(M))($locally defined on $U)$ such that
\begin{equation*}
\nabla T_\alpha=\omega_\alpha^\beta T_\beta,\nabla T_{\bar{\alpha}}=\omega_{\bar{\alpha}}^{\bar{\beta}} T_{\bar{\beta}},\nabla T=0.
\end{equation*}
The $\omega_\alpha^\beta$ satisfies the following equations:
\begin{equation*}
d\theta^\beta=\theta^\alpha\wedge \omega_\alpha^\beta+\theta\wedge\tau^\beta,
\end{equation*}
\begin{equation*}
0=\tau_\alpha\wedge \theta^\alpha,
\end{equation*}
\begin{equation*}
0=\omega_\alpha^\beta+\omega_{\bar{\beta}}^{\bar{\alpha}},
\end{equation*}
where $\tau_\alpha=A_{\alpha\gamma}\theta^\gamma$ with $A_{\alpha\gamma}=A_{\gamma\alpha}$. The curvature of Tanaka-Webster connection, expressed in terms of the coframe $\{\theta,\theta^\alpha,\theta^{\bar{\alpha}}\}$, is
\begin{equation*}
\Pi_\beta^\alpha=d\omega_\beta^\alpha-\omega_\beta^\gamma\wedge \omega_\gamma^\alpha.
\end{equation*}
Webster showed that $\Pi_\beta^\alpha$ can be written
\begin{equation*}
\Pi_\beta^\alpha=R_{\beta\rho\bar{\sigma}}^\alpha\theta^\rho\wedge\theta^{\bar{\sigma}}+W_{\beta\rho}^\alpha\theta^\rho\wedge\theta-W_{\beta\bar{\rho}}^\alpha\theta^{\bar{\rho}}\wedge\theta
+i \theta_\beta\wedge\tau^\alpha-i\tau_\beta\wedge \theta^\alpha.
\end{equation*}
We will denote components of covariant derivatives with indices preceded by a comma; thus write $A_{\alpha\beta,\gamma}$. The indices $\{0,\alpha,\bar{\alpha}\}$ indicate derivatives with respect to $\{T,T_\alpha,T_{\bar{\alpha}}\}$. For derivatives of a scalar function, we will often omit the comma, for instance, $u_\alpha=T_\alpha(u),u_{\alpha\bar{\beta}}=T_\alpha(u_{\bar{\beta}})-\omega_\alpha^\gamma(T_{\bar{\beta}})u_\gamma$, $u_{\alpha\beta\gamma}=
T_\alpha(u_{\beta\gamma})-\omega_\alpha^\lambda(T_{\beta})u_{\lambda\gamma}-u_{\beta\mu}\omega_\alpha^\mu(T_{\gamma})$. Notice that the notations we used are somewhat different from the authors in \cite{CKLT}.

For a real function $u$, the subgradient $\nabla^\mathrm{b}$ is defined by $\nabla^\mathrm{b} u\in \Gamma(H(M))$ and $\langle Z,\nabla^\mathrm{b} u\rangle=du(Z)$ for all $Z\in \Gamma(H(M))$. Locally $\nabla^\mathrm{b} u=u_\alpha T_{\bar{\alpha}}+u_{\bar{\alpha}}T_\alpha$. We can use the connection to define the subhessian by
\begin{equation*}
\pi_\mathrm{b} \nabla^2u=u_{\alpha\beta}\theta^\alpha\otimes\theta^\beta+u_{\alpha\bar{\beta}}\theta^\alpha\otimes\theta^{\bar{\beta}}
+u_{\bar{\alpha}\beta}\theta^{\bar{\alpha}}\otimes\theta^\beta+u_{\bar{\alpha}\bar{\beta}}\theta^{\bar{\alpha}}\otimes\theta^{\bar{\beta}}.
\end{equation*}
In particular, $|\nabla^\mathrm{b} u|^2=2u_\alpha u_{\bar{\alpha}}$, $|\pi_\mathrm{b} \nabla^2 u|^2=2(u_{\alpha\beta}u_{\bar{\alpha}\bar{\beta}}+u_{\alpha\bar{\beta}}u_{\bar{\alpha}\beta})$. Also, we define the sub-Laplacian by
\begin{equation*}
\Delta_\mathrm{b}u=tr(\pi_\mathrm{b} \nabla^2u)=\sum_\alpha(u_{\alpha\bar{\alpha}}+u_{\bar{\alpha}\alpha}).
\end{equation*}
For simplicity, we always denote $T(\varphi)$ by $\varphi_0$ for any smooth function $\varphi$ on $M$.

Next we recall the following commutation relations. Let us set $\{T_A\}=\{T, T_\alpha, T_{\bar{\alpha}}\}$, where $A \in\{0,1, \ldots, n, \overline{1}, \ldots, \bar{n}\}$ and $T_0=T)$. Let $\varphi$ be a smooth function, then we have 
\begin{align}
	& \varphi_{\alpha \beta}=\varphi_{\beta \alpha}, \notag\\
	& \varphi_{\alpha \bar{\beta}}=\varphi_{\bar{\beta} \alpha}- i h_{\alpha \bar{\beta}} f_0, \notag\\
	& \varphi_{\alpha 0}=\varphi_{0 \alpha}+A_\alpha^{\bar{\beta}} \varphi_{\bar{\beta}}.\label{P.1}
\end{align}
and 
\begin{align*}
& \varphi_{\lambda \mu \alpha}=\varphi_{\mu \lambda \alpha}+i(\varphi_\mu A_{\alpha \lambda}-\varphi_\lambda A_{\alpha \mu}), \\
& \varphi_{\lambda \bar{\mu} \alpha}=\varphi_{\bar{\mu} \lambda \alpha}- i \varphi_{0 \alpha} h_{\lambda \bar{\mu}}-\varphi_\beta R_\alpha{ }^\beta \lambda \bar{\mu}, \\
& \varphi_{\bar{\lambda} \bar{\mu} \alpha}=\varphi_{\bar{\mu} \bar{\lambda} \alpha}+ i \varphi_\beta(h_{\alpha \bar{\mu}} A_{\bar{\lambda}}^\beta-h_{\alpha \bar{\lambda}} A_{\bar{\mu}}^\beta), \\
& \varphi_{\mu 0 \alpha}=\varphi_{0 \mu \alpha}+f_{\bar{\beta} \alpha} A_\mu^{\bar{\beta}}-\varphi_\beta W_{\alpha \mu}^\beta, \\
& \varphi_{\bar{\mu} 0 \alpha}=\varphi_{0 \bar{\mu} \alpha}+\varphi_{\beta \alpha} A_{\bar{\mu}}^\beta+\varphi_\beta W_{\alpha \bar{\mu}}^\beta.
\end{align*}
\begin{definition}\label{CRL}
We say that $(M,J,\theta)$ satisfies the CR sub-Lapacian comparison property if there exist constants $C_0=C_0(n)\geq 0,\tilde{k}\geq 0$ such that	
$$
\Delta_\mathrm{b} r \leq C_0(1/r+\sqrt{\tilde{k}}),
$$
in the distribution sense.
\end{definition}
\section{Proof of the subparabolic Li-Yau type estimate in theorem $\ref{thm}$}
This section is devoted to the proof of the nonlinear version of the Li-Yau estimate in Theorem $\ref{thm}$. As the proof is quite involved and requires several intermediate steps, for the sake of clarity and convenience, we break this into three subsection, focusing first on deriving and establishing some of the necessary tools and identities and then finalizing the proof in the last subsection.
\subsection{CR Bochner formula and some basic lemmas}
We first recall the following CR version of Bochner formula from Greenleaf \cite{G} in a pseudo-Hermitian $(2n+1)$-manifold.
\begin{lemma}\label{3Lem1}
For a smooth real-valued function $\varphi$ on $M$,
\begin{align*}
\Delta_\mathrm{b} |\nabla^\mathrm{b} \varphi|^2 &= 2|\pi_\mathrm{b}\nabla^2 \varphi|^2+2\langle\nabla^\mathrm{b} \varphi,\nabla^\mathrm{b} \Delta_\mathrm{b}\varphi\rangle \\
& \quad+2(2\operatorname{Ric}-(n-2)\operatorname{Tor})(\nabla^{1,0}\varphi,\nabla^{1,0}\varphi)+4\langle J\nabla^\mathrm{b}\varphi,\nabla^\mathrm{b} \varphi_0\rangle.
\end{align*}
\end{lemma}
From Lemma $\ref{3Lem1}$, the authors in \cite{CKLT} obtained the following CR Bochner-type estimate.
\begin{lemma}\label{lemBochesti}
For a smooth real-valued function $\varphi$ on $M$ and any $\nu>0$, we have
\begin{align}
\Delta_\mathrm{b} |\nabla^\mathrm{b} \varphi|^2 & \geq 4\bigg(\sum_{\alpha,\beta}|\varphi_{\alpha\beta}|^2+\sum_{\alpha,\beta=1,\alpha\neq \beta}^n|\varphi_{\alpha\bar{\beta}}|^2\bigg)+\frac{1}{n}(\Delta_\mathrm{b}\varphi)^2+n\varphi_0^2  \notag\\
&\quad+2\langle\nabla^\mathrm{b} \varphi,\nabla^\mathrm{b} \Delta_\mathrm{b}\varphi\rangle +2(2\operatorname{Ric}-(n-2)\operatorname{Tor})(\nabla^{1,0}\varphi,\nabla^{1,0}\varphi)\notag\\
&\quad-\frac{2}{\nu}|\nabla^\mathrm{b}\varphi|^2-2\nu|\nabla^\mathrm{b} \varphi_0|^2.\label{Bochnerestimate}
\end{align}
\end{lemma}
The following lemma gives a relationship between $\Delta_\mathrm{b} \varphi_0$ and $(\Delta_\mathrm{b}\varphi)_0$ for a smooth function $\varphi$.
\begin{lemma}[\cite{CKLT}]
For a smooth real-valued function $\varphi$ on $M$, we have
\begin{equation}\label{31.1}
\Delta_\mathrm{b} \varphi_0=(\Delta_\mathrm{b}\varphi)_0+2\sum_{\alpha,\beta=1}^n[(\varphi_{\bar{\beta}}A_{\alpha\beta})_{\bar{\alpha}}+(\varphi_\beta A_{\bar{\alpha}\bar{\beta}})_\alpha].
\end{equation}
\end{lemma}
Now, we define
\begin{equation*}
V(\varphi)=\sum_{\alpha,\beta=1}^n[(A_{\alpha\beta}\varphi_{\bar{\beta}})_{\bar{\alpha}}+(A_{\bar{\alpha}\bar{\beta}}\varphi_\beta)_\alpha
+A_{\alpha\beta}\varphi_{\bar{\beta}}\varphi_{\bar{\alpha}}+A_{\bar{\alpha}\bar{\beta}}\varphi_\beta\varphi_\alpha].
\end{equation*}
Let $u$ be a positive smooth solution of $(\ref{equation})$ on $M\times [0,+\infty)$ and denote
\begin{equation*}
f=\log u.
\end{equation*}
It is easy to check that $f$ satisfies the equation
\begin{equation}\label{31.2}
(\Delta_\mathrm{b}-\partial_t)f=-af-|\nabla^\mathrm{b} f|^2.
\end{equation}
\begin{lemma}\label{3Lem2}
Let $u$ be a positive smooth solution of $(\ref{equation})$ on $M\times [0,+\infty)$ with $f=\log u$. Then
\begin{equation*}
\Delta_\mathrm{b} f_0=f_{0t}-af_0-2\langle \nabla^\mathrm{b} f_0,\nabla^\mathrm{b} f\rangle+2V(f),
\end{equation*}
where $f_{0t}=T(f_t)$.
\begin{proof}
Using $(\ref{31.1})$, $(\ref{31.2})$ adnd the commutation relation $(\ref{P.1})$, we have
\begin{align*}
\Delta_\mathrm{b} f_0 &=(f_t-af-|\nabla^\mathrm{b}f|^2)_0+ 2\sum_{\alpha,\beta=1}^n[(\varphi_{\bar{\beta}}A_{\alpha\beta})_{\bar{\alpha}}+(\varphi_\beta A_{\bar{\alpha}\bar{\beta}})_\alpha]\\
&=f_{0t}-af_0-2\langle \nabla^\mathrm{b} f_0,\nabla^\mathrm{b} f \rangle+2V(f).
\end{align*}
\end{proof}
\end{lemma}
\begin{lemma}[\cite{CKLT}]\label{3Lem3}
Let $u$ be a positive smooth function on $M\times [0,+\infty)$ with $f=\log u$. Suppose that
\begin{equation*}
[\Delta_\mathrm{b},T]u=0.
\end{equation*}
Then
\begin{equation*}
V(f)=0.
\end{equation*}
\end{lemma}
From Lemmas $\ref{3Lem2}$ and $\ref{3Lem3}$, we easily obtain the following conclusion.
\begin{corollary}\label{3Cor1}
Let $u$ be a positive smooth solution of $(\ref{equation})$ on $M\times [0,+\infty)$ with $f=\log u$. Suppose that
\begin{equation*}
[\Delta_\mathrm{b},T]u=0.
\end{equation*}
Then $f$ satisfies
\begin{equation}\label{31.3}
\Delta_\mathrm{b} f_0=f_{0t}-af_0-2\langle \nabla^\mathrm{b} f_0,\nabla^\mathrm{b} f\rangle.
\end{equation}
\end{corollary}

\subsection{A subparabolic inequality}
In this subsection we introduce a Harnack quantity built out of the solution $u$ and establish a subparabolic inequality under the operator $\partial_t - \Delta_\mathrm{b}$.
\begin{lemma}
Let $u$ be a positive smooth solution of $(\ref{equation})$ on $B_{2R}(x_0)\times [0,+\infty)$ and let $F=F(x,t)$ be defined by
\begin{equation}\label{32.2}
F(x,t)=t(|\nabla^\mathrm{b}f|^2(x,t)+\alpha a f(x,t)-\alpha f_t(x,t)+\beta tf_0^2(x,t)),\;\; t\geq 0
\end{equation} 
where $f=\log u$ and $\alpha>1,\beta>0$ are constants. Suppose that
\begin{equation*}
[\Delta_\mathrm{b},T]u=0.
\end{equation*}
Then $F$ satisfies
\begin{align*}
\Delta_\mathrm{b}F-F_t&=-2\langle\nabla^\mathrm{b} f,\nabla^\mathrm{b} F\rangle-\frac{F}{t}-aF +t\bigg[2|\pi_\mathrm{b}\nabla^2 f|^2+4\langle J\nabla^\mathrm{b}f,\nabla^\mathrm{b} f_0\rangle\notag\\
&\quad+2(2\operatorname{Ric}-(n-2)\operatorname{Tor})(\nabla^{1,0}f,\nabla^{1,0}f)+2\beta t |\nabla^\mathrm{b} f_0|^2\notag\\
&\quad +(\alpha-1)a|\nabla^\mathrm{b} f|^2-\beta taf_0^2-\beta f_0^2\bigg],
\end{align*}
on $B_{2R}(x_0)\times (0,+\infty)$.
\begin{proof}
Referring to the equation for $u$, an easy calculation shows that $f$ in turn satisfies the equation 
\begin{equation}\label{32.4}
 (\Delta_\mathrm{b}-\partial_t)f=-af-|\nabla^\mathrm{b} f|^2.
\end{equation} 
Moreover, using $(\ref{32.2})$ and $(\ref{32.4})$, it is easily seen that the following relation emerges between $F$, $|\nabla^\mathrm{b} f|^2$ and $\Delta_\mathrm{b}f$:
\begin{align}
\Delta_\mathrm{b}f&=f_t-af-\bigg(\frac{F}{t}-\alpha a f+\alpha f_t-\beta tf_0^2\bigg)\notag\\
&=-\frac{F}{t}+(\alpha-1)af-(\alpha-1)f_t+\beta t f_0^2.\label{32.5}
\end{align}

Now having these identities and relations in place we next proceed onto applying the operator $\partial_t - \Delta_\mathrm{b}$ to the Harnack quantity $F$ given by $(\ref{32.2})$. Towards this end, we first note that 
\begin{equation}\label{32.6}
\Delta_\mathrm{b}F=t(\Delta_\mathrm{b}|\nabla^\mathrm{b}f|^2+\alpha a \Delta_\mathrm{b}f-\alpha \Delta_\mathrm{b}f_t+2\beta tf_0\Delta_\mathrm{b}(f_0)+2\beta t |\nabla^\mathrm{b} f_0|^2).
\end{equation}
As for the first term on the right by recalling the CR Bochner formula as applied to $f$, we have
\begin{align*}
\Delta_\mathrm{b} |\nabla^\mathrm{b} f|^2 &= 2|\pi_\mathrm{b}\nabla^2 f|^2+2\langle\nabla^\mathrm{b} f,\nabla^\mathrm{b} \Delta_\mathrm{b}f\rangle \\
& \quad+2(2\operatorname{Ric}-(n-2)\operatorname{Tor})(\nabla^{1,0}f,\nabla^{1,0}f)+4\langle J\nabla^\mathrm{b}f,\nabla^\mathrm{b} f_0\rangle,
\end{align*}
and so upon substituting back in $(\ref{32.6})$ and making use of $(31.3)$ this gives
\begin{align}
\Delta_\mathrm{b}F&=t\bigg[2|\pi_\mathrm{b}\nabla^2 f|^2+2\langle\nabla^\mathrm{b} f,\nabla^\mathrm{b} \Delta_\mathrm{b}f\rangle+2(2\operatorname{Ric}-(n-2)\operatorname{Tor})(\nabla^{1,0}f,\nabla^{1,0}f)\notag\\
&\quad\quad+4\langle J\nabla^\mathrm{b}f,\nabla^\mathrm{b} f_0\rangle+\alpha a \Delta_\mathrm{b}f-\alpha \Delta_\mathrm{b}f_t+2\beta t |\nabla^\mathrm{b} f_0|^2\notag\\
&\quad\quad +2\beta tf_0f_{0t}-2\beta taf_0^2-4\beta tf_0\langle \nabla^\mathrm{b} f_0,\nabla^\mathrm{b} f\rangle\bigg].\label{32.7}
\end{align}

Now referring to the sum on the right the contributions of the second, fifth and sixth terms, modulo a factor $t$ and upon using $(\ref{32.4})$ and $($(\ref{32.5})$)$ can be simplified and re-written as,
\begin{align*}
2\langle\nabla^\mathrm{b} f,\nabla^\mathrm{b} \Delta_\mathrm{b}f\rangle&+\alpha a \Delta_\mathrm{b}f-\alpha \Delta_\mathrm{b}f_t\\
&=2\bigg \langle\nabla^\mathrm{b} f,\nabla^\mathrm{b} \bigg(-\frac{F}{t}+(\alpha-1)af-(\alpha-1)f_t+\beta t f_0^2\bigg)\bigg\rangle\\
&\quad +\alpha a \bigg(-\frac{F}{t}+(\alpha-1)af-(\alpha-1)f_t+\beta t f_0^2\bigg) \\
&\quad -\alpha f_{tt}+a\alpha f_t+2\alpha \langle \nabla^\mathrm{b} f,\nabla^\mathrm{b} f_t\rangle)\\
&=-\frac{2}{t}\langle\nabla^\mathrm{b} f,\nabla^\mathrm{b} F\rangle+2(\alpha-1)a|\nabla^\mathrm{b} f|^2+2\langle\nabla^\mathrm{b} f,\nabla^\mathrm{b} f_t\rangle\\
&\quad +4\beta tf_0 \langle\nabla^\mathrm{b} f,\nabla^\mathrm{b}f_0\rangle-\frac{\alpha a F}{t}+\alpha a^2 (\alpha-1)f\\
&\quad -(\alpha-2)\alpha af_t+\alpha a\beta tf_0^2-\alpha f_{tt}.
\end{align*}
Therefore substituting this back into $(\ref{32.7})$ and rearranging the equality give
\begin{align}
\Delta_\mathrm{b}F&=-2\langle\nabla^\mathrm{b} f,\nabla^\mathrm{b} F\rangle +t\bigg[2|\pi_\mathrm{b}\nabla^2 f|^2+2(2\operatorname{Ric}-(n-2)\operatorname{Tor})(\nabla^{1,0}f,\nabla^{1,0}f)\notag\\
&\quad+4\langle J\nabla^\mathrm{b}f,\nabla^\mathrm{b} f_0\rangle+2\beta t |\nabla^\mathrm{b} f_0|^2+2(\alpha-1)a|\nabla^\mathrm{b} f|^2+2\langle\nabla^\mathrm{b} f,\nabla^\mathrm{b} f_t\rangle\notag\\
&\quad -\frac{\alpha a F}{t}+\alpha a^2 (\alpha-1)f+2\beta tf_0f_{0t}+(\alpha-2)(\beta taf_0^2-\alpha af_t)-\alpha f_{tt}\bigg].\label{32.8}
\end{align}
Combining $(\ref{32.8})$ and
\begin{align*}
F_t &= \frac{F}{t}+t(|\nabla^\mathrm{b}f|^2+\alpha a f-\alpha f_t+\beta tf_0^2)_t\notag\\
&=\frac{F}{t}+t[2\langle\nabla^\mathrm{b}f,\nabla^\mathrm{b} f_t\rangle+\alpha af_t-\alpha f_{tt}+\beta f_0^2+2\beta tf_0 f_{0t}],
\end{align*}
we obtain
\begin{align*}
\Delta_\mathrm{b}F-F_t&=-2\langle\nabla^\mathrm{b} f,\nabla^\mathrm{b} F\rangle-\frac{F}{t} +t\bigg[2|\pi_\mathrm{b}\nabla^2 f|^2+4\langle J\nabla^\mathrm{b}f,\nabla^\mathrm{b} f_0\rangle\notag\\
&\quad+2(2\operatorname{Ric}-(n-2)\operatorname{Tor})(\nabla^{1,0}f,\nabla^{1,0}f)+2\beta t |\nabla^\mathrm{b} f_0|^2\notag\\
&\quad +2(\alpha-1)a|\nabla^\mathrm{b} f|^2-\frac{\alpha a F}{t}+\alpha a^2 (\alpha-1)f\\
&\quad+(\alpha-2)\beta taf_0^2-(\alpha-1)\alpha af_t-\beta f_0^2\bigg]\\
&=-2\langle\nabla^\mathrm{b} f,\nabla^\mathrm{b} F\rangle-\frac{F}{t} +t\bigg[2|\pi_\mathrm{b}\nabla^2 f|^2+4\langle J\nabla^\mathrm{b}f,\nabla^\mathrm{b} f_0\rangle\notag\\
&\quad+2(2\operatorname{Ric}-(n-2)\operatorname{Tor})(\nabla^{1,0}f,\nabla^{1,0}f)+2\beta t |\nabla^\mathrm{b} f_0|^2\notag\\
&\quad +(\alpha-1)a|\nabla^\mathrm{b} f|^2-\frac{\alpha a F}{t}+(\alpha-2)\beta taf_0^2\\
&\quad+(\alpha-1)a\bigg(|\nabla^\mathrm{b} f|^2+\alpha a f-\alpha f_t\bigg)-\beta f_0^2\bigg].
\end{align*}
Using $(\ref{32.2})$ again and rearranging the equality leads to the desired conclusion.
\end{proof}
\end{lemma}
\begin{lemma}\label{3Lem4}
Let $(M^{2n+1},J,\theta)$ be a complete pseudo-Hermitian manifold with
\begin{equation}\label{32.9}
(2\operatorname{Ric}-(n-2)\operatorname{Tor})(Z,Z)\geq-2k|Z|^2,\text{ for all }Z\in T_{1,0}(M), \text{ on }B_{2R}(x_0),
\end{equation}
where $k\geq0$. Let $u$ be a positive smooth solution of $(\ref{equation})$ on $B_{2R}(x_0)\times [0,+\infty)$ and let $F=F(x,t)$ be as in $(\ref{32.2})$. Suppose that
\begin{equation*}
[\Delta_\mathrm{b},T]u=0.
\end{equation*}
Then
\begin{align}
\Delta_\mathrm{b} F-F_t & \geq -2\langle \nabla^\mathrm{b} f,\nabla^\mathrm{b} F\rangle-\frac{F}{t}-aF \notag\\
& \quad+t\bigg[\frac{1}{n}(\Delta_\mathrm{b} f)^2+(n-\beta-ta\beta)f_0^2\notag\\
& \quad+\bigg((\alpha-1)a-2k-\frac{2}{\beta t}\bigg)|\nabla^\mathrm{b} f|^2\bigg],\label{32.10}
\end{align}
on $B_{2R}(x_0)\times (0,+\infty)$.
\begin{proof}
The conclusion follows at once by recalling $(\ref{Bochnerestimate})$, making note of $(\ref{32.9})$ and taking $\nu=\beta t$.
\end{proof}
\end{lemma}
\subsection{Proof of the local estimate in Theorem $\ref{thm}$}
Having all the ingredients and necessary tools at our disposal we now come to the proof of the main estimate. The idea of the proof is to combine the inequality established in Lemma $\ref{3Lem4}$ along with a localization argument by using a suitable cut-off function. The estimates on the cut-off function in turn makes use of the CR sub-Laplacian comparison property as will be described in detail in the course of the proof. We pick a reference point $x_0\in M$ and fix $R>0,\tau>0$. We denote by $r(x)=d(x,x_0)$ the geodesic radial variable in reference to $x_0$. For the sake of localization we consider first a function $\chi=\chi(r)$ on the half-line $r\geq 0$ $($see Lemma $\ref{3Lem5}$ below$)$ and then for $x\in M$ set
\begin{equation}\label{33.1}
\phi(x)=\chi\bigg(\frac{r(x)}{R}\bigg).
\end{equation}
The existence of $\chi$ as used in $(\ref{33.1})$ and its properties is granted by the following straightforward and statement.
\begin{lemma}\label{3Lem5}
There exists a function $\chi:[0,+\infty)\rightarrow \mathbb{R}$ verifying the following properties:\\
$i)$ $\chi$ is of class $C^2[0,+\infty)$;\\
$ii)$ $0\leq \chi(r)\leq $ for $0\leq r<\infty$ with $\chi=1$ for $r\leq 1$ and $\chi=0$ for $s\geq 2$;\\
$iii)$ $\chi$ is non-increasing and additionally, for suitable constants $\epsilon, \upsilon>0$, satisfies the bounds
\begin{equation}\label{33.2}
-\epsilon\chi^{1/2}(r)\leq \chi^\prime(r),\;\;\; \text{and} \;\;\;\chi^{\prime\prime}\geq -\upsilon,
\end{equation}
on the half-line $[0,+\infty)$.
\end{lemma}
\noindent \textbf{Case 1: $a\leq0$.}

It is evident from ii) that $\phi=1$ for when $0\leq r(x)\leq R$ and $\phi=0$ for when $r(x)\geq 2R$. Let us now consider the spatially localized function $\phi F$ where $F$ is as in $(\ref{32.2})$. We denote by $(z,s)$ the point where this function attains its maximum over the compact set $\bar{B}_{2R}(x_0)\times [0,\tau]$. By virtue of Calabi's standard argument\cite{C} we can assume that $z$ is not in the cut locus of $x_0$ and so $\phi$ is smooth at $z$ for the application of the maximum principle. Additionally, we can assume that $[\phi F](z,s)>0$ as otherwise the desired estimate is trivially true as a result of $F\leq 0$. It thus follows that $s>0$ and $z\in B_{2R}(x_0)$ and so at the maximum point $(z,s)$ we have the relations 
\begin{equation}\label{33.3}
\nabla^\mathrm{b} (\phi F)=0,\Delta_\mathrm{b}(\phi F)\leq 0,F_t \geq 0.
\end{equation}
Starting with the basic identity $\Delta_\mathrm{b}(\phi F)=\phi \Delta_\mathrm{b} F+2\langle \nabla^\mathrm{b} \phi,\nabla^\mathrm{b} F\rangle + F\Delta_\mathrm{b} \phi$ and making note of the relations $(\ref{33.3})$ at maximum point $(z,s)$ we can write
\begin{align}
0&\geq \Delta_\mathrm{b}(\phi F)=\phi \Delta_\mathrm{b} F+2\phi^{-1}\langle \nabla^\mathrm{b} \phi,\nabla^\mathrm{b} (\phi F)\rangle -2(|\nabla^\mathrm{b}\phi|^2/\phi)F+ F\Delta_\mathrm{b} \phi\notag\\
&\geq \phi \Delta_\mathrm{b} F -2(|\nabla^\mathrm{b}\phi|^2/\phi)F+ F\Delta_\mathrm{b} \phi.\label{33.4}
\end{align}
Now from $(\ref{33.1})$ we deduce $\nabla^\mathrm{b}\phi=(\chi^\prime/R)\nabla^\mathrm{b}r$, $\Delta_\mathrm{b} \phi=\chi^{\prime\prime}|\nabla^\mathrm{b}r|^2/R^2+\chi^\prime \Delta_\mathrm{b} r /R$ and so 
\begin{equation}\label{33.5}
|\nabla^\mathrm{b}\phi|^2/\phi \leq \frac{\epsilon^2}{R^2},\;\;\Delta_\mathrm{b} \phi\geq -\frac{\upsilon}{R^2}-\frac{\epsilon C_0}{R^2}(1+\sqrt{\tilde{k}}R),
\end{equation}
where we use $(\ref{33.2})$ and the CR sub-Laplacian comparison property $(\ref{CRL})$. Thus returning to $(\ref{33.4})$, invoking $(\ref{32.10})$ and making note of $(\ref{33.3})$, $(\ref{33.5})$, we obtain, at the maximum point $(z,s)$, the inequality 
\begin{align}
0&\geq \phi \Delta_\mathrm{b} F -2(|\nabla^\mathrm{b}\phi|^2/\phi)F+ F\Delta_\mathrm{b} \phi\notag\\
&\geq  -\bigg(\frac{\upsilon}{R^2}+\frac{\epsilon C_0}{R^2}(1+\sqrt{\tilde{k}}R)\bigg) F-2\frac{\epsilon^2}{R^2}F\notag\\
&\quad +\phi\bigg\{-2\langle \nabla^\mathrm{b} f,\nabla^\mathrm{b} F\rangle-\frac{F}{s}-aF\notag\\
&\quad \quad +s\bigg[\frac{1}{n}(\Delta_\mathrm{b} f)^2+(n-\beta-sa\beta)f_0^2\notag\\
&\quad \quad +\bigg((\alpha-1)a-2k-\frac{2}{\beta s}\bigg)|\nabla^\mathrm{b} f|^2\bigg]\bigg\}.\label{33.6}
\end{align}
We now bound the first term on the second line in the last inequality, using $(\ref{33.3})$ and $(\ref{33.5})$, we have,
\begin{equation*}
2\phi \langle \nabla^\mathrm{b} f,\nabla^\mathrm{b} F\rangle=-2\langle \nabla^\mathrm{b} f,\nabla^\mathrm{b} \phi\rangle F\leq 2 |\nabla^\mathrm{b}\phi||\nabla^\mathrm{b}f|F\leq \frac{2\epsilon}{R}\phi^{1/2}|\nabla^\mathrm{b}f|F.
\end{equation*} 
For simplify, set 
$$B=\frac{\upsilon}{R^2}+\frac{\epsilon C_0}{R^2}(1+\sqrt{\tilde{k}}R)+2\frac{\epsilon^2}{R^2},\text{ and } \mu=|\nabla^\mathrm{b} f|^2(z,s)/F(z,s)\geq 0.$$
As a result using the above in $(\ref{33.6})$ and rearranging terms we have
\begin{align}
 & BF +\frac{2\epsilon F^{3/2}}{R}\phi^{1/2}\mu^{1/2}+\phi F/s+a\phi F\notag\\
&\geq s\phi\bigg[\frac{1}{n}(\Delta_\mathrm{b} f)^2+(n-\beta-sa\beta)f_0^2 +\bigg((\alpha-1)a-2k-\frac{2}{\beta s}\bigg)F\mu\bigg].\label{33.7}
\end{align}
at $(z,s)$. From the definition of $F,\mu$ and the equation $(\ref{32.4})$, we have, at $(z,s)$
\begin{align}
(\Delta_\mathrm{b} f)^2&=(f_t-af -|\nabla^\mathrm{b}f|^2)\notag\\
&=\bigg(\frac{F}{\alpha s}+(1-\alpha^{-1})|\nabla^\mathrm{b}f|^2-\alpha^{-1}\beta sf_0^2\bigg)^2\notag\\
&\geq \frac{F^2(1+(\alpha-1)\mu s)^2}{\alpha^2 s^2}-2\frac{[1+(\alpha-1)\mu s]F\beta  f_0^2}{\alpha^2 }.\label{33.8}
\end{align}
Multiplying both sides of $(\ref{33.7})$ by $s\phi$ and using $(\ref{33.8})$ and the fact that $0\leq\phi\leq1$ and $a\leq0$, we obtain at $(z,s)$,
\begin{align}
&s\phi BF+2\frac{\epsilon F^{3/2}}{R}\mu^{1/2}s\phi^{3/2}+F\phi\notag\\
& \geq s^2\phi^2\bigg[\frac{1}{n}(\Delta_\mathrm{b} f)^2+(n-\beta)f_0^2-\bigg((\alpha-1)|a|+2k+\frac{2}{\beta s}\bigg)F\mu\bigg]\notag\\
&\geq\frac{\phi^2F^2}{n}\frac{(1+(\alpha-1)\mu s)^2}{\alpha^2 }-2\frac{s^2\phi^2}{n}\beta f_0^2F\frac{1+(\alpha-1)\mu s}{\alpha^2}\notag\\
&\quad+s^2\phi^2\bigg[(n-\beta)f_0^2-\bigg((\alpha-1)|a|+2k+\frac{2}{\beta s}\bigg)F\mu\bigg]\notag\\
&=\frac{\phi^2F^2(1+(\alpha-1)\mu s)^2}{n\alpha^2 }-s^2\phi^2\bigg((\alpha-1)|a|+2k+\frac{2}{\beta s}\bigg)F\mu\notag\\
&\quad+\bigg[(n-\beta)-2\frac{\beta F}{n}\frac{1+(\alpha-1)\mu s}{\alpha^2}\bigg]f_0^2s^2\phi^2.\label{33.9}
\end{align}
Now, set
\begin{equation*}
\beta=\frac{n}{1+\displaystyle{\frac{2 F(z,s)}{n}\frac{1+(\alpha-1)\mu s}{\alpha^2}}}\text{ and }\lambda=\phi(z)F(z,s).
\end{equation*}
Then we have at $(z,s)$
\begin{align}
sB\lambda+2\frac{\epsilon \lambda^{3/2}}{R}\mu^{1/2}s+\lambda
&\geq\frac{\lambda^2(1+(\alpha-1)\mu s)^2}{n\alpha^2 }-s^2\phi^2((\alpha-1)|a|+2k)F\mu\notag\\
&\quad -\frac{ 2 s\lambda\mu}{n}-\frac{4 s\lambda^2\mu}{n^2}\frac{1+(\alpha-1)\mu s}{\alpha^2}\notag\\
&\geq \frac{\lambda^2(1+(\alpha-1)\mu s)[(1+(\alpha-1)\mu s)-4 s\mu/n]}{n\alpha^2 }\notag\\
&\quad -s^2\lambda\mu((\alpha-1)|a|+2k)-\frac{ 2 s\lambda\mu}{n}.\label{33.10}
\end{align}
For simplicity, set
\begin{equation*}
A=(1+(\alpha-1)\mu s)[(1+(\alpha-1)\mu s)-4 s\mu/n].
\end{equation*}
By the assumption $(\alpha-1)n\geq4$, we have $A\geq 1$. For any fixed $0<\delta<1$, at $(z,s)$
\begin{equation}\label{33.11}
2\frac{\epsilon \lambda^{3/2}}{R}\mu^{1/2}s\leq (1-\delta)\frac{\lambda^2}{n\alpha^2}A+\frac{\lambda\epsilon^2s^2\mu \alpha^2 n}{(1-\delta)A R^2}.
\end{equation}
Using $(\ref{33.10})$ and $(\ref{33.11})$, we have
\begin{equation}\label{33.12}
\frac{\delta\lambda}{n\alpha^2}\leq \frac{s^2\mu((\alpha-1)|a|+2k)+sB+1}{A}+\frac{ 2 s\mu}{nA}+\frac{\epsilon^2s^2\mu \alpha^2 n}{(1-\delta)A^2 R^2}.
\end{equation}
We now proceed onto bounding the full expression on the right-hand side in the last inequality. Towards this end dealing with the first term, 
\begin{align}
\frac{s^2\mu((\alpha-1)|a|+2k)+sB+1}{(1+(\alpha-1)\mu s)[(1+(\alpha-1)\mu s)-4 s\mu/n]} & \leq sB+\frac{s^2\mu((\alpha-1)|a|+2k)+1}{1+(\alpha-1)\mu s}  \notag\\
&\leq sB+1+\frac{s((\alpha-1)|a|+2k)}{\alpha-1}.\label{33.13}
\end{align}
In much the same way regarding the subsequent terms, we have
\begin{equation}\label{33.14}
\frac{ 2 s\mu}{n(1+(\alpha-1)\mu s)[(1+(\alpha-1)\mu s)-4 s\mu/n]}\leq \frac{2}{n(\alpha-1)},
\end{equation}
and
\begin{equation}\label{33.15}
\frac{\epsilon^2s^2\mu \alpha^2 n}{(1-\delta)(1+(\alpha-1)\mu s)^2[(1+(\alpha-1)\mu s)-4 s\mu/n]^2 R^2}\leq\frac{\epsilon^2s \alpha^2 n}{2(1-\delta)(\alpha-1)R^2}.
\end{equation}

Now referring to $(\ref{33.12})$, making use of the bounds obtained in $(\ref{33.13})$-$(\ref{33.15})$, we obtain
\begin{equation*}
\lambda\leq \frac{ n\alpha^2}{\delta}\bigg[sB+1+\frac{s((\alpha-1)|a|+2k)}{\alpha-1}+\frac{2}{n(\alpha-1)}+\frac{\epsilon^2s \alpha^2 n}{2(1-\delta)(\alpha-1)R^2}\bigg].
\end{equation*}
Recalling that the definition of $\lambda$ and $B$, we have
\begin{align*}
\sup_{B_R(x_0)} &|\nabla^\mathrm{b}f|^2(x,\tau)+\alpha a f(x,\tau)-\alpha f_t(x,\tau)\leq
\frac{n\alpha^2}{\delta\tau}\bigg(1+\frac{2}{n(\alpha-1)}\bigg)\\
&\quad+ \frac{n\alpha^2}{\delta}\bigg[\bigg(\frac{\upsilon}{R^2}+\frac{\epsilon C_0}{R^2}(1+\sqrt{\tilde{k}}R)+\frac{2\epsilon^2}{R^2}\bigg)+\frac{((\alpha-1)|a|+2k)}{\alpha-1}\\
&\quad+\frac{\epsilon^2 \alpha^2 n}{2(1-\delta)(\alpha-1)R^2}\bigg],
\end{align*}
We complete the proof of $(i)$ in Theorem $\ref{thm}$.

\noindent \textbf{Case 2: $a>0$.}

We adopt the same notation as above. We will also assume $\phi F$ achieves its positive maximum at $(z,s)\in B_{2R}(x_0)\times [0,\tau]$. Using a similar means of deriving $(\ref{33.9})$, we obtain the following estimate at point $(z,s)$:
\begin{align*}
&s\phi BF+2\frac{\epsilon F^{3/2}}{R}\mu^{1/2}s\phi^{3/2}+F\phi+asF\phi\\
& \geq\frac{\phi^2F^2(1+(\alpha-1)\mu s)^2}{n\alpha^2 }-s^2\phi^2\bigg(2k+\frac{2}{\beta s}\bigg)F\mu\\
&\quad+\bigg[(n-\beta-sa\beta)-2\frac{\beta F}{n}\frac{1+(\alpha-1)\mu s}{\alpha^2}\bigg]f_0^2s^2\phi^2.
\end{align*}
Now, set
\begin{equation*}
\beta=\frac{n}{1+sa+\displaystyle{\frac{2 F(z,s)}{n}\frac{1+(\alpha-1)\mu s}{\alpha^2}}}\text{ and }\lambda=\phi(z)F(z,s).
\end{equation*}
Note that $a>0$. Then we have at $(z,s)$
\begin{align}
sB\lambda+2\frac{\epsilon \lambda^{3/2}}{R}\mu^{1/2}s&+\lambda+as\lambda
\geq\frac{\lambda^2(1+(\alpha-1)\mu s)^2}{n\alpha^2 }-2s^2k\lambda\mu\notag\\
&\quad -\frac{ 2 s\lambda\mu(1+sa)}{n}-\frac{4 s\lambda^2\mu}{n^2}\frac{1+(\alpha-1)\mu s}{\alpha^2}\notag\\
&=\frac{\lambda^2(1+(\alpha-1)\mu s)[(1+(\alpha-1)\mu s)-4 s\mu/n]}{n\alpha^2 }\notag\\
&\quad -2s^2k\lambda\mu-\frac{ 2 s\lambda\mu(1+sa)}{n}.\label{33.16}
\end{align}
By $(\ref{33.11})$, the inequality $(\ref{33.16})$ becomes
\begin{equation}\label{33.17}
\frac{\delta\lambda}{n\alpha^2 } \leq\frac{ 2s^2k\mu+1+as+sB}{A}+\frac{ 2 s\mu(1+sa)}{nA}+\frac{\epsilon^2s^2\mu \alpha^2 n}{(1-\delta)A^2 R^2}.
\end{equation}
Similarlly, we now proceed onto bounding the full expression on the right-hand side in the last inequality: 
\begin{align}
\frac{s^2\mu((\alpha-1)|a|+2k)+1}{(1+(\alpha-1)\mu s)[(1+(\alpha-1)\mu s)-4 s\mu/n]} & \leq \frac{s^2\mu((\alpha-1)|a|+2k)+1}{1+(\alpha-1)\mu s} \notag \\
&\leq 1+\frac{s((\alpha-1)|a|+2k)}{\alpha-1},\label{33.18}
\end{align}
\begin{equation}\label{33.19}
\frac{ 2 s\mu}{n(1+(\alpha-1)\mu s)[(1+(\alpha-1)\mu s)-4 s\mu/n]}\leq \frac{2}{n(\alpha-1)},
\end{equation}
and
\begin{equation}\label{33.20}
\frac{\epsilon^2s^2\mu \alpha^2 n}{(1-\delta)(1+(\alpha-1)\mu s)^2[(1+(\alpha-1)\mu s)-4 s\mu/n]^2 R^2}\leq\frac{\epsilon^2s \alpha^2 n}{2(1-\delta)(\alpha-1)R^2}.
\end{equation}

Now referring to $(\ref{33.12})$, making use of the bounds obtained in $(\ref{33.19})$-$(\ref{33.20})$, we obtain
\begin{equation*}
\lambda\leq \frac{ n\alpha^2}{\delta}\bigg[sB+1+\frac{s((\alpha-1)|a|+2k)}{\alpha-1}+\frac{2}{n(\alpha-1)}+\frac{\epsilon^2s \alpha^2 n}{2(1-\delta)(\alpha-1)R^2}\bigg].
\end{equation*}
Recalling that the definition of $\lambda$ and $B$, we have
\begin{align*}
\sup_{B_R(x_0)} &|\nabla^\mathrm{b}f|^2(x,\tau)+\alpha a f(x,\tau)-\alpha f_t(x,\tau)\leq
\frac{n\alpha^2}{\delta \tau}\bigg(1+\frac{2}{n(\alpha-1)}\bigg)\\
&\quad+ \frac{n\alpha^2}{\delta}\bigg[\bigg(\frac{\upsilon}{R^2}+\frac{\epsilon C_0}{R^2}(1+\sqrt{\tilde{k}}R)+\frac{2\epsilon^2}{R^2}\bigg)+\frac{((\alpha-1)|a|+2k)}{\alpha-1}\\
&\quad+\frac{\epsilon^2 \alpha^2 n}{2(1-\delta)(\alpha-1)R^2}\bigg].
\end{align*}
We complete the proof of $(ii)$ in Theorem $\ref{thm}$.
\section{Proof of corollary $\ref{cor}$}
Assume that $u(x)$ is a positive smooth solution to the equation
\begin{equation*}
\Delta_\mathrm{b}u+au\log u=0,\;\text{ on }M,
\end{equation*}
satisfying 
\begin{equation*}
[\Delta_\mathrm{b},T]u=0.
\end{equation*}

\noindent	\textbf{Case 1.} $a<0$.

Since $k=0$ and $u(x)$ is independent of $t$, by $(i)$ of Theorem $\ref{thm}$, we have on $B_R(x_0)$,
\begin{equation}\label{4.1}
|\nabla^\mathrm{b}f|^2+\alpha a f\leq \frac{n\alpha^2}{\delta}\bigg[\bigg(\frac{\upsilon}{R^2}+\frac{\epsilon C_0}{R^2}(1+\sqrt{\tilde{k}}R)+\frac{2\epsilon^2}{R^2}\bigg)+|a|+
\frac{\epsilon^2 \alpha^2 n}{2(1-\delta)(\alpha-1)R^2}\bigg],
\end{equation}
where $\alpha\geq 4/n+1$ and $0<\delta<1$ are two real constants, and $\epsilon>0$ and $\upsilon>0$ are some fixed constants. Letting $R\rightarrow +\infty$, we have by $(\ref{4.1})$,
\begin{equation}\label{4.2}
\alpha a \log u(x)\leq -\frac{an\alpha^2}{\delta}.
\end{equation}
Since $a<0$, taking $\alpha=4/n+1$ and $\delta\rightarrow 1$, we have by $(\ref{4.2})$,
\begin{equation*}
u(x)\geq e^{-4-n}\text{ for all }x\in M.
\end{equation*}

\noindent	\textbf{Case 2.} $a>0$.

Using the same method of deriving $(\ref{4.2})$, we have by $(ii)$ of Theorem $\ref{thm}$,
\begin{equation}\label{4.3}
\alpha a \log u(x)\leq \frac{n\alpha^2}{\delta}\bigg(a+\frac{2a}{n(\alpha-1)}\bigg).
\end{equation}
Since $a>0$, taking $\alpha=4/n+1$ and $\delta\rightarrow 1$, we have $(\ref{4.3})$,
\begin{equation*}
u(x)\leq e^{3(4+n)/2}\text{ for all }x\in M.
\end{equation*}

\section{Proof of the parabolic Harnack inequality in Theorem $\ref{thm2}$}
With the aid of the estimates established in Theorem $\ref{thm}$ we can now prove the desired subparabolic Harnack inequality in Theorem $\ref{thm2}$. Towards this end it suffices to integrate the former estimate along suitable space-times curves in $B_{R/2}(x_0)\subset M\times (0,+\infty)$. Let us first move on to rewriting the inequality $(\ref{0.1})$ and $(\ref{0.2})$ as follows.
\begin{align}
\frac{|\nabla^\mathrm{b}u|^2}{\alpha u^2}&+ a \log u-\frac{u_t}{u}\leq
\frac{n\alpha}{\delta t}\bigg(1+\frac{2}{n(\alpha-1)}\bigg)\notag\\
&\quad+ \frac{n\alpha}{\delta}\bigg[\bigg(\frac{\upsilon}{R^2}+\frac{\epsilon C_0}{R^2}(1+\sqrt{\tilde{k}}R)+\frac{2\epsilon^2}{R^2}\bigg)+|a|+\frac{(2a^++2nk)}{n(\alpha-1)}\notag\\
&\quad+
\frac{\epsilon^2 \alpha^2 n}{2(1-\delta)(\alpha-1)R^2}\bigg].\label{5.1}
\end{align}
Let $S=S(n,\alpha,\delta,\upsilon,R,\epsilon,C_0,k,\tilde{k},a)$ be defined by 
\begin{align}
S=&-\frac{n\alpha}{\delta}\bigg[\bigg(\frac{\upsilon}{R^2}+\frac{\epsilon C_0}{R^2}(1+\sqrt{\tilde{k}}R)+\frac{2\epsilon^2}{R^2}\bigg)+|a|+\frac{(2a^++2nk)}{n(\alpha-1)}\notag\\
&\quad+
\frac{\epsilon^2 \alpha^2 n}{2(1-\delta)(\alpha-1)R^2}\bigg]+a q_1,\label{5.2}
\end{align}
where $q_1=\sup_{B_R(x_0)\times [0,\tau]}\log u$. It follows from $(\ref{5.1})$ that 
$$u_t/u\geq \frac{|\nabla^\mathrm{b}u|^2}{\alpha u^2}-\frac{n\alpha}{\delta t}\bigg(1+\frac{2}{n(\alpha-1)}\bigg)+S.$$ Suppose that $\gamma: [t_1,t_2]\rightarrow M$ is a lengthy curve from $x_1$ to $x_2$. Recall that $\gamma$ is a lengthy curve if $\dot{\gamma}(t)\in H(M)_{\gamma(t)}$ for every $t$. By a theorem of W.L. Chow \cite{Chow}, any two points in $M$ can be joined by a lengthy curve. Using the above it is seen that 
\begin{align*}
\frac{d}{dt}\log u(\gamma(t),t)&=\bigg\langle \frac{\nabla^\mathrm{b}u}{u},\dot{\gamma}(t)\bigg\rangle+\frac{u_t}{u}\\
&\geq\bigg\langle \frac{\nabla^\mathrm{b}u}{u},\dot{\gamma}(t)\bigg\rangle+\frac{|\nabla^\mathrm{b}u|^2}{\alpha u^2}-\frac{n\alpha}{\delta t}\bigg(1+\frac{2}{n(\alpha-1)}\bigg)+S\\
&=\alpha^{-1}\bigg|\frac{\nabla^\mathrm{b}u}{u}+\alpha \frac{\dot{\gamma}(t)}{2}\bigg|^2-\alpha \frac{|\dot{\gamma}(t)|^2}{4}-\frac{n\alpha}{\delta t}\bigg(1+\frac{2}{n(\alpha-1)}\bigg)+S\\
&\geq -\alpha \frac{|\dot{\gamma}(t)|^2}{4}-\frac{n\alpha}{\delta t}\bigg(1+\frac{2}{n(\alpha-1)}\bigg)+S.
\end{align*}
Therefore integrating the above inequality gives 
\begin{align*}
\log \frac{u(x_2,t_2)}{u(x_1,t_1)}&=\int_{t_1}^{t_2}\frac{d}{dt}\log u(\gamma(t),t) dt\\
&\geq \int_{t_1}^{t_2}-\alpha \frac{|\dot{\gamma}(t)|^2}{4}-\int_{t_1}^{t_2}\frac{n\alpha}{\delta t}\bigg(1+\frac{2}{n(\alpha-1)}\bigg) dt+S(t_2-t_1)\\
&=-\frac{n\alpha}{\delta }\bigg(1+\frac{2}{n(\alpha-1)}\bigg)\log(t_2/t_1)-\frac{\alpha}{4} \int_{t_1}^{t_2}|\dot{\gamma}(t)|^2 dt+S(t_2-t_1).
\end{align*}
Hence upon exponentiating we have
\begin{equation*}
\frac{u(x_2,t_2)}{u(x_1,t_1)}\geq \bigg(\frac{t_2}{t_1}\bigg)^{-\frac{\alpha(n\alpha-n+2)}{\delta(\alpha-1) }}\exp\bigg[-\frac{\alpha}{4} \int_{t_1}^{t_2}|\dot{\gamma}(t)|^2 dt\bigg]\exp[S(t_2-t_1)],
\end{equation*}
or upon rearranging terms and rescaling the integral:
\begin{equation*}
u(x_2,t_2)\geq u(x_1,t_1)\bigg(\frac{t_2}{t_1}\bigg)^{-\frac{\alpha(n\alpha-n+2)}{\delta(\alpha-1) }}e^{-\alpha L(x_1,x_2,t_2-t_1)}e^{S(t_2-t_1)}
\end{equation*}
where 
$$ L(x_1,x_2,t_2-t_1)=\inf_{\gamma}\bigg[\frac{1}{4(t_2-t_1)}\int_0^1 |\dot{\gamma}(t)|^2 dt\bigg].$$

\bibliographystyle{amsplain}

\end{document}